\documentclass[11pt]{article}
\usepackage{amssymb,amsmath,latexsym}
\usepackage{color}

\usepackage{epsfig}

\newtheorem{thm}{Theorem}[section]
\newtheorem{defn}[thm]{Definition}

\newtheorem{lemma}[thm]{Lemma}

\newtheorem{remark}[thm]{Remark}
\newtheorem{example}[thm]{Example}
\newtheorem{assumption}[thm]{Assumption}

\newcommand{\pf}{\noindent{\bf Proof.} }
\def\qed{{\hfill $\Box$ \bigskip}}

\def\loc{{\rm loc}}

\newcommand\cbrk{\text{$]$\kern-.15em$]$}}
\newcommand\opar{\text{\,\raise.2ex\hbox{${\scriptstyle
|}$}\kern-.34em$($}}
\newcommand\cpar{\text{$)$\kern-.34em\raise.2ex\hbox{${\scriptstyle |}$}}\,}

\def\wh{\widehat}
\def\wt{\widetilde}
\def\<{\langle}
\def\>{\rangle}
\def\eps{\varepsilon}
\def\E{{\mathbb E}}

\newcommand\bL{\mathbb{L}}
\newcommand\bR{\mathbb{R}}
\newcommand\bH{\mathbb{H}}
\newcommand\bZ{\mathbb{Z}}

\newcommand\bD{\mathbb{D}}
\newcommand\bS{\mathbb{S}}

\newcommand\cB{\mathcal{B}}

\newcommand\cF{\mathcal{F}}
\newcommand\cH{\mathcal{H}}

\newcommand\cP{\mathcal{P}}
\newcommand\cR{\mathcal{R}}

\def\wh{\widehat}
\def\wt{\widetilde}

\newcommand{\mysection}[1]{\section{#1}
\setcounter{equation}{0}}

\begin{document}

\title{\bf An $L^2$-theory on
 SPDE driven by L\'e{}vy processes}

\date{}

\author{Zhen-Qing Chen\footnote{Department of Mathematics, University of Washington,
Seattle, WA 98195, USA, \,\, zchen@math.washington.edu. The research
of this author is supported in part by NSF Grant DMS-0906743.}
\qquad \hbox{\rm and} \qquad Kyeong-Hun Kim\footnote{Department of
Mathematics, Korea University, 1 Anam-dong, Sungbuk-gu, Seoul, South
Korea 136-701, \,\, kyeonghun@korea.ac.kr. The research of this
author is supported by the Korean Research Foundation Grant funded
by the Korean Government 20090087117}}

\maketitle

\begin{abstract}
In this paper  we develop an $L_2$-theory for stochastic  partial
differential equations driven by L\'e{}vy processes.
 The coefficients of the equations are random functions depending on time
and space variables, and no smoothness assumption of the
coefficients is assumed.

\vspace*{.125in}


\vspace*{.125in}


\end{abstract}

\mysection{Introduction}

 In this article
we study the $L^2$-theory of
  stochastic partial differential
equations of the following type:
\begin{equation} \label{e:1.1}
du= \left(\frac{\partial}{\partial x_i} \left( a^{ij}u_{x^j}+
\bar{b}^iu \right)+b^iu_{x^i}+  cu+ f \right) dt + \left(
\sigma^{ik} u_{x^i}+\mu^{k}u+g^{k} \right) dZ^{k}_{t}
\end{equation}
given for $t\geq0$  and $x\in \bR^d$. Here $\{Z^k_t$, $k=2, 1,
\cdots\}$ are independent one-dimensional L\'evy processes, $i$ and
$j$ go from $1$ to $d$ with the summation convention on $i,j,k$
being enforced.
%new
For example, the second term in the right hand side of \eqref{e:1.1}
should be understood as
$$  \sum_{k\geq 1} \left( \sum_{i=1}^d
   a^{ik}u_{x^i}+\mu^k u +g^k\right) dZ^k_t.
$$
%end new
The coefficients $a^{ij}$, $b^i$, $c, \sigma^{ik},
\mu^{k}$ and the free terms $f, g^k$ are {\bf{random}} functions
depending on $(t,x)$.

Stochastic partial differential equations (SPDEs) of type \eqref{e:1.1}
arise naturally in applications when the objects are subject to randomness
and high variability.
The purpose of this
paper is to investigate the existence and uniqueness of pathwise
solutions
 to \eqref{e:1.1} and to study the regularity of the solutions.

If $Z^k_t$ are independent one-dimensional Wiener processes, then
general $L^p$-theory of the equation   has been well studied. An
$L^p$-theory of SPDEs with Wiener processes defined on $\bR^n$ was
first introduced by Krylov in \cite{Kr99}, and in \cite{KL1} and
\cite{KL2} Krylov and Lototsky developed an $L^p$-theory of such
equations with constant coefficients defined on half space
$\bR^n_+$. Later in many articles (see \cite{KK}, \cite{Kim03} and
references therein) these results were extended for SPDEs with
variable coefficients defined on bounded domains of $\bR^n$.

However very little is known when $Z^k_t$ are general discontinuous
L\'evy processes. In \cite{CZ}, existence and uniqueness of  weak
(or martingale) solutions as well as pathwise solutions to the
following SPDE
\begin{equation}  \label{e:1.2}
 du={\cal A}u dt +\sum_{k=1}^n g^k(u)dZ^k_t,
 \end{equation}
driven by L\'evy processes is studied,
where ${\cal A}$ is the generator of certain semigroup on a Hilbert space
% and the functions $g^k$, $k=1, \cdots, n$,
%  are non-random and satisfy certain continuity condition.
$H$  and $g^k$, $k=1, \cdots, n$,
are non-random maps from $H$ to $H$ that satisfy certain continuity condition.

See the Introduction of \cite{CZ} for a brief discussion on other
related work SPDE driven by Poisson random measure or stable noises,
including \cite{AWZ, F, M, RZ}.
%Note that unlike in our article equation \eqref{e:1.2} has only
%{\bf{non-random}} coefficients, which is   {\bf{independent}} of $t$.
 Note that  maps $g^k$, $k=1,\cdots, n$, in \eqref{e:1.2}
are {\bf{non-random}} coefficients and are {\bf{independent}} of $t$,
while $g^k$'s in \eqref{e:1.1} to be considered in this paper
are random and time dependent but are given a priori that
do not depend on solution u.
 Moreover
 %the first derivatives of solutions do not appear  in
 no derivatives of the solution $u$ appear in
the stochastic part of equation (\ref{e:1.2}).

Our main result, Theorem \ref{main thm},
 is presented
and proved in
section \ref{section 2}. Here we
%prove that if
show that if each $Z^k_t$ has finite second moment, i.e., if
\begin{equation}
                 \label{eqn 222}
\int_\bR z^2 \nu_k(dz) <\infty \quad \hbox{for every } k\geq 1,
\end{equation}
where $\nu^k$ is the L\'evy measure of $Z^k$, then equation
(\ref{e:1.1}) admits a unique solution in
$\bH^1(T):=L^2(\Omega\times [0,T],W^1_2)$ and the $\bH^1(T)$-norm of
the solution is controlled by the $L^2$-norm of $f$ and $g$. In
section \ref{section extension} we give two extensions of Theorem
\ref{main thm}. First we
 develop
an $L^2$-theory for a certain type of nonlinear equations.
%And second,
  Second,
we weaken condition (\ref{eqn 222}) by assuming that it
holds only for sufficiently large $k$ (thus it can be dropped if
only finitely many processes $Z^k$ appear in the equation) and prove
that the equation has unique pathwise $W^1_2$-valued solution.

As usual, throughout this paper,
$\bR^{d}$ stands for the Euclidean space of points $x=(x^{1},...,x^{d})$.
 For $i=1,...,d$, multi-indices $\alpha=(\alpha_{1},...,\alpha_{d})$,
$\alpha_{i}\in\{0,1,2,...\}$, and functions $u(x)$,  we set
$$
u_{x^{i}}=\partial u/\partial x^{i}=D_{i}u,\quad
D^{\alpha}u=D_{1}^{\alpha_{1}}\cdot...\cdot D^{\alpha_{d}}_{d}u,
\quad|\alpha|=\alpha_{1}+...+\alpha_{d}.
$$
We also use the notation $D^m$ for a partial derivative of order $m$
with respect to $x$. If we write $c=c(...)$,
%this
 it
means that the
constant $c$ depends only on what are in parenthesis.

\mysection{Main results}
                               \label{section 2}
Let $(\Omega,\cF,P)$ be a complete probability space equipped with a
 filtration $(\cF_{t},t\geq0)$ satisfying the usual condition.
   We assume that
on $\Omega$ we are given independent one-dimensional Levy processes
 $Z^{1}_{t},Z^{2}_{t},...$
relative to $\{\cF_{t},t\geq0\}$. Let $\cP$ be the predictable
$\sigma$-field generated by $\{\cF_{t},t\geq0\}$.

For $t\geq 0$ and $A\in \cB(\bR \setminus \{0\})$,  define
$$
N_k(t,A)=\left\{0\leq s\leq t; \, Z^k_s-Z^k_{s-} \in A \right\}, \quad \wt
{N}_k(t,A)=N_k(t,A)-t\nu_k(A)
$$
where $\nu_k(A):=\E [N_k(1,A)]$ is the L\'e{}vy measure of $Z^k$. By
L\'e{}vy-It\^o decomposition,  there exist constants $\alpha^k,
\beta^k$ and Brownian motion $B^k$ so that
\begin{equation}
                           \label{eqn 01.15.10}
Z^k_t=\alpha^kt +\beta^k B^k_t+\int_{|z|<1}z \wt
{N}_k(t,dz)+\int_{|z|\geq 1} z N_k(t, dz).
\end{equation}

\begin{assumption}        \label{A2.1}
{\rm (i)}  For each $k\geq 1$,
\begin{equation}
                 \label{eqn e.1}
\wh {c}_{k}:= \left[\int_\bR z^2 \nu_k(dz)\right]^{1/2} <\infty.
\end{equation}

{\rm (ii)} There exist constants $\delta,K >0$ so that for every
$t>0$, $x\in \bR^d$ and
 $\omega \in \Omega$,
\begin{equation}    \label{e:2.2}
\delta|\xi|^2 \leq ( a^{ij}-\alpha^{ij} )\xi^i\xi^j\leq
a^{ij}\xi^i\xi^j\leq K|\xi|^2, \quad \forall \xi\in \bR^d,
\end{equation}
where $\alpha^{ij}:= \frac{1}{2}
%summation convention is inforce so removed: \sum_{k=1}^\infty
(\wh {c}_{k}^2+\beta^2_k) \sigma^{ik}\sigma^{jk}$.
%moved from below.
 Here $i$ and $j$ go from $1$ to $d$, and $k$
runs through $\{1,2,\cdots\}$.
\end{assumption}

Recall that throughout the article,  summation convention is used.
%end moving
 Due to (\ref{eqn e.1}), $\int_{|z|>1} |z| N_k(1, dz)<\infty$, and thus by absorbing $\wt
\alpha_k:=\int_{|z|>1} z N_k(1, dz)$ into $\alpha_k$ we can rewrite
(\ref{eqn 01.15.10}) as
$$
Z^k_t=\tilde{\alpha}_k t +\beta_k B^k_t+\int_{\bR^1}z \wt
{N}_k(t,dz).
$$

%moved below
%Without loss of generality, we assume that $\tilde{\alpha}^k=0$.
 For $d\geq 1$, consider the equation for random function
 $u(t, x)$ on $\Omega \times [0, \infty)\times \bR^d$:
\begin{equation}      \label{5.15.0}
du= \left(\frac{\partial}{\partial x_i} \left( a^{ij}u_{x^j}+
\bar{b}^iu \right)+b^iu_{x^i}+  cu+ f \right) dt + \left(
\sigma^{ik} u_{x^i}+\mu^{k}u+g^{k} \right) dZ^{k}_{t}
\end{equation}
in the weak sense. See Definition \ref{D:2.5} below.
% Here    $i$ and $j$ go from $1$ to $d$, and $k$
%runs through $\{1,2,\cdots\}$. Throughout the article such summation
%convention is used.
The coefficients $a^{ij},\bar{b}^i$, $b^i$, $c$,
$\sigma^{ik}$, $\mu^k$ and the free terms $f,g^k$ are random
functions depending on $t>0$ and $x\in \bR^d$.
%Note that if $\tilde{\alpha}_k\neq 0$, then it is enough to
Without loss of generality, we assume that $\tilde{\alpha}^k=0$,
since otherwise we can simply
 move the term
$\sum_{k}\tilde{\alpha}_k\left(\sigma^{ik}u_{x^i} +\nu^k u+
g^k\right)\,dt $ from the stochastic part to the deterministic one.

\begin{remark}\rm
 Conditions (\ref{eqn e.1}) and (\ref{e:2.2}) will be weakened in section \ref{section
extension}. In particular, one can completely drop the condition
(\ref{eqn e.1}) if there are only finitely many processes $Z^k_t$ in
equation (\ref{5.15.0}).
\end{remark}

For $n=0,1,2,...$, let
$$
H^n:=\left\{u\in L^2(\bR^d):\,  Du,...,D^n u\in L^2 (\bR^d) \right\},
$$
which is equipped with norm $\| u\|_{H^n}:= \left( \sum_{k=0}^n \|D^k u\|_{L^2(\bR^d)}^2 \right)^{1/2}$.
Here $Du := (\frac{\partial u}{\partial x_1}, \cdots, \frac{\partial u}{\partial x_d})$
denotes the gradient of $u$ in the distributional sense, $D^2u$ denotes the collection
of all second derivatives of $u$ in the distribution sense, and so on.
Let $H^{-n}:=(H^n)^{*}$ be its topological dual,  and
$\cP^{dP\times dt}$ be the completion of $\cP$ with respect to  $dP
\times dt$. For $n\in \bZ$ and $T>0$, we write $u\in \bH^{n}(T)$ if
$u$ is an $H^{n}$-valued
 $\cP^{dP\times dt}$-measurable  process defined on
$\Omega \times [0,T]$ so that
$$
\|u\|_{\bH^n(T)}:= \left( \E \left[\int^T_0\|u (t, \cdot )\|^2_{H^{n}}\,dt \right] \right)^{1/2}<\infty.
$$
Denote $\bL(T):=\bH^0(T)$. For an $\ell^2$-valued processes
$g=(g_1,g_2,...)$, we say $g\in \bL(T,\ell^2)$ if  $g^k\in \bL(T)$
for every $k\geq 1$ and
$$
\|g\|_{\bL(T,\ell^2)}:=
% \left( \sum_{k=1}^{\infty}(\beta^2_k+\wh{c}_k^2)
 \sum_{k=1}^{\infty}(\beta^2_k+\wh{c}_k^2) \left(
 \E \left[\int^T_0\|g^k\|^2_{L^2}\,dt\right] \right)^{1/2} <\infty.
$$
Finally
we use $U_2$ to denote the family of  $L^2(\bR^d)$-valued $\cF_0$-measurable random variables $u_0$ having
$$
\|u_0\|_{U_2}:=\left( \E  \left[\|u_0\|^2_{L^2}
\right]\right)^{1/2}<\infty.
$$

\begin{remark}\label{R:2.4}\rm \begin{description}
%\item{(i)} Let $\cM^n(T)$ denote the set of all $H^n$-valued $\cF \otimes \cB([0,T])$-measurable, $\{\cF_t\}_{t\geq 0}$-adapted processes $u(t)$ on $[0, T]$ so that
%$$
%\E \left[\int^T_0\|u (t) \|^2_{H^n}dt \right]<\infty.
%$$
%By Theorem 2.8.2 in \cite{kr}, $ \cM^n(T) \subset \bH^n(T)$.
%ZC question: is the above an equality?

\item{(i)} Since we assume $\tilde{\alpha}_k=0$,
 $Z^k  $ is a square integrable martingale, whose quadratic
variational process will be denoted as $[Z^k]$. By L\'evy system,
the predictable dual projection $\< Z^k\>$ of $[Z^k]$ is given by
$\<Z^k\>_t = (\wh c_k^2+\beta^2_k)\, t$.
 For every
process $H$ in $L^2(\Omega\times [0,T])$, which has a predictable
%version, $M_t:=\int_0^t H_s dZ^k_s$ is again a square integrable
 $d P\times dt$-version $\wt H$, $M_t=\int^t_0 H_s dZ^k_s:=\int_0^t \wt H_s dZ^k_s$
 is well defined and is independent of the choice of such $\wt H$,
 and $M$ is a
martingale with
$$ \E \left[ M_t^2 \right] = \E \left[ \int_0^t H_s^2 \, d[Z^k]_s\right]
 =  (\beta^2_k+\wh c_k^2) \,  \E \left[ \int_0^t H_s^2 ds\right], \quad t\leq T.
 $$
%new
 We will simply denote $M$ by $\int_0^\cdot H_s dZ_s^k$.
 %end new
%  So for any
 For
 $g=(g_1, g_2, \cdots )\in \bL(T,\ell_2)$ and $\phi\in
C_c^\infty(\bR^d)$,
$$
\sum_{k=1}^{\infty}\int^T_0(\beta^2_k+\wh{c}_k^2)(g^k,\phi)^2ds \leq
\|\phi\|^2_{L^2}\|g\|^2_{\bL(T,\ell^2)}  <\infty \quad \hbox{a.s.}
$$
Therefore the series of stochastic integral $\sum_{k=1}^{\infty}
\int_0^t (g^k,\phi)\, dZ^k_t$ converges uniformly in
%$t$ in probability on $[0,T]$.
$t\in [0, T]$ in probability.

\item{(ii)} In many other articles, the equation of the type
$$
du=(Au+f)dt+g(u(t-))dZ_t
$$
has been studied. The expression $u(t-)$ is used so that it is
predictable and the integral $\int^t_0 g(u(t-))dZ_t$ becomes a
martingale. Such notation is not used  in (\ref{5.15.0}),  because
%as mentioned in (ii)
 by (i) and (ii),
stochastic integral can be defined for a
process $H$ in $L_2(\Omega\times [0,T])$ which has a predictable
version $\tilde{H}$, and
$$
\int^t_0 H(s)dZ_s = \int^t_0 \tilde{H}(s)dZ_s.
$$

\end{description}
\end{remark}

\begin{defn}\label{D:2.5}
We say $u \in \cH^1(T)$ if $u\in \bH^1(T)$, $u$  is right continuous
having left limits  in $L^2$ $a.s.$ with $u(0)\in U_2$,
 and for some $f\in
\bH^{-1}(T)$ and $ g=(g_1, g_2, \cdots )\in \bL(T,\ell^2)$
$$
du(t)=f(t) dt +  g^k(t) dZ^k_t \qquad \hbox{for  } 0\leq t\leq T
$$
 in the sense of distributions; that is, for any $\phi\in C^\infty_c (\bR^d)$, the equality
\begin{equation}
                \label{eqn 11.16}
(u(t),\phi)=(u(0),\phi)+\int^t_0(f(s),\phi)ds +
\sum_{k=1}^\infty\int^t_0(g^k (s),\phi)dZ^k_s
\end{equation}
holds for all $t\leq T$ $a.s.$. In this case, we write
$$
\bD u:=f, \quad \bS u:=g,
$$
and define
$$
\|u\|_{\cH^1(T)}:=\|u\|_{\bH^1(T)}+\|\bD u\|_{\bH^{-1}(T)}+\|\bS
u\|_{\bL(T,\ell^2)}+\|u(0)\|_{U_2}.
$$
\end{defn}

\begin{lemma}       \label{lemma 11.16}
Let $u\in \cH^1(T)$, then

{\rm (i)} for any $\phi\in H^1$, $(u(t),\phi)$ is progressively
measurable, right continuous
% with left limits ;
 having left limits ;

{\rm (ii)} for each fixed $t>0$,  $u(t)=u(t-)$ in $L^2$ a.s.

\end{lemma}
\pf (i) follows immediately from (\ref{eqn 11.16}).

(ii). By assumption $u(t-)$ exists. Let $\{\phi_n, :\phi_n \in H^1, n=1,2,...\}$ be a orthonormal
basis in $L^2 (\bR^d)$. Then the process $t\mapsto (u(t-),\phi_n)$
is predictable by (i). Since $\int^t_0(g^k,\phi_n)dZ^{k}_t$ is
stochastically continuous, we have for each fixed $t$ and $n\geq 1$,
$(u(t),\phi_n)=(u(t-),\phi_n)$ a.s.   Therefore
$$
u(t-)=\sum_n(u(t-),\phi_n)\phi_n=u(t) \quad \hbox{a.s.}
$$
 The lemma is now proved.
\qed

\begin{thm}
                        \label{T:2.7}
The space $\cH^1(T)$ is a Banach space and
\begin{equation}
                                         \label{e:2.4}
\E \left[ \sup_{t\leq T}\|u(t)\|^2_{L^2} \right] \leq
c\left(
%\|u_x\|^2_{\bL(T)}
 \|D u \|^2_{\bL(T)}
+\|\bD u\|^2_{\bH^{-1}(T)}+\|\bS
u\|^2_{\bL(T,\ell^2)}+ \E \|u(0)\|^2_{L^2}\right),
\end{equation}
where $c$ is independent of $u$ and $T$.

\end{thm}

\pf First we prove \eqref{e:2.4}. Let $u(0)=u_0$ and $du=f dt+
 g^k dZ^k_t$. Then for any $\phi\in C^{\infty}_0$,
\begin{equation}
             \label{eqn 5.19.1}
     (u(t),\phi)=(u(0),\phi)+\int^t_0(f(s),\phi)ds+\int^t_0(g^k(s),\phi)dZ^k_t
     \end{equation}
     for all $t\leq T$ (a.s.).
     %Write $f\in \bH^{-1}(T)$ as
 For $f\in \bH^{-1}(T)$, we can write it as
%$$
%f=f_0+\sum_{i=1}^d \frac{\partial}{\partial x_i} f_i \quad\hbox{with
%}  f_i\in \bL(T)
%$$
%so that
%$$
%\sum_{i=0}^d\|f_i\|_{\bL(T)}\leq c\|f\|_{\bH^{-1}(T)}.
%$$
$$
f=f_0+\sum_{i=1}^d \frac{\partial}{\partial x_i} f_i
$$
with $f_i\in \bL(T)$ for $0\leq i\leq d$ and
$$
\sum_{i=0}^d\|f_i\|_{\bL(T)}\leq c\|f\|_{\bH^{-1}(T)}.
$$

Indeed, since $f=(1-\Delta)(1-\Delta)^{-1}f$ and
$(1-\Delta)^{-1}:H^n \to H^{n+2}$ is an isometry, we can take
$$
f_0=(1-\Delta)^{-1}f \qquad\hbox{and} \qquad f_i=-\frac{\partial
f_0}{\partial x^i} \quad \hbox{for } \ i=1,2,...,d.
$$
Take a nonnegative function $\psi\in C^{\infty}_0(B_1(0))$ with unit
integral, and for $\varepsilon>0$ define
$\psi_{\varepsilon}(x)=\varepsilon^{-d}\psi(x/\varepsilon)$.
 For any generalized function $u$, define
$u^{(\varepsilon)}(x)=u*
%\zeta_{\varepsilon}(x)
 \psi_\varepsilon (x)
:=(u(\cdot),\psi_{\varepsilon}(x-\cdot))$,
then $u^{(\varepsilon)}(x)$ is infinitely differentiable function of
$x$.  By plugging $\psi_{\varepsilon}(x-\cdot)$ instead of $\phi$ in
(\ref{eqn 5.19.1}),
$$
u^{(\varepsilon)}(t,x)=u^{(\varepsilon)}(0,x)+\int^t_0
(f^{(\varepsilon)}_0+D_if^{(\varepsilon)}_i)dt+\int^t_0g^{(\varepsilon)k}dZ^k_t.
$$

By  taking $\varepsilon \to 0$, one can easily show  that
(\ref{e:2.4}) holds true if for any $\varepsilon>0$ it holds with
$u^{(\varepsilon)}, u_0^{(\varepsilon)}, f^{(\varepsilon)},
g^{(\varepsilon)}$ in place of $u,u_0,f,g$, respectively. Thus we
may assume that $u,f,g$ are infinitely differentiable in $x$, and
therefore (a.s.)
\begin{equation}
                  \label{eqn 1}
u(t)=u_0+\int^t_0fdt+\int^t_0g^kdZ^k_t, \quad \forall t\leq T.
\end{equation}
The stochastic integral in (\ref{eqn 1}) doesn't change if we
replace $g$ by its predictable version, thus we also assume that $g$
is predictable.

Applying Ito's formula to $|u(t)|^2$ (cf. \cite{HWY}) and
integrating over $\bR^d$, we have
 \begin{eqnarray}\label{e:2.5}
\|u(t)\|^2_{L^2}&=&
 \| u_0\|_2^2 + 2 \int_0^t (u(s), f(s)) ds + \sum_k
\beta^2_k\int^t_0|g^k(s)|^2_{L^2}ds \nonumber\\
&& \ \
 +2\sum_k\int^t_0(u(s-),g^k (s))dZ^k_s +\sum_k\sum_{0<s\leq
t}\|g^k(s)\Delta Z^k_s\|^2_{L^2} \nonumber \\
&=& \|u_0\|_{L^2}^2+2\int^t_0 \left( (u(s), \, f_0 (s))
  -\sum_{i=1}^d(u_{x^i}(s),f_i(s))\right)ds
  %new
  + \sum_k
\beta^2_k\int^t_0|g^k(s)|^2_{L^2}ds
%end new
 \nonumber\\
&& \ \
 +2\sum_k\int^t_0(u(s-),g^k (s))dZ^k_s +\sum_k\sum_{0<s\leq
t}\|g^k(s)\Delta Z^k_s\|^2_{L^2},
\end{eqnarray}
where we have used the fact
% that
 that $Z^k$'s are independent and so
 with probability one at most one of
the $Z^1_s, Z^2_s\cdots $ can jump at any given time. By virtue of
the L\'e{}vy system of the L\'e{}vy process $Z^k_s$, it follows that
\begin{equation}\label{e:3.6}
 \sum_{0<s\leq t}\|g^k(s)\Delta Z^k_s\|^2_{L^2}=M_t^k+ \wh {c}_{k}^2 \int_0^t
\|g^k\|^2_{L^2} ds,
\end{equation}
where $M^k$ is a purely discontinuous square integrable martingale
with
$$
 M^k_t - M^k_{t-}=\| g^k(t)\Delta Z^k_t\|^2_{L^2} \qquad \hbox{for } t> 0 .
$$
%As easy to check,
It is easy to see that for every $\eps>0$, there is a constant
$c(\eps)>0$, independent of $u$ and $f_i$'s such that
\begin{eqnarray*}
 && \E \left[\sup_{t\leq T} \Big|\int^t_0 \big( (u(s), f_0
(s))-\sum_{i=1}^d(u_{x^i}(s),
%f^i(s)
 f_i(x)
)\big)ds \Big|\right] \\
% &\leq& \|D u\|^2_{\bL(T)}
%  +\varepsilon \E\sup_{t\leq T}\|u(t)\|^2_{L^2}+c(\varepsilon)\sum_{i=0}^d\|f^i\|^2_{\bL(T)} \\
%  &\leq& \|D u\|^2_{\bL(T)}
%  +\varepsilon \E\sup_{t\leq T}\|u(t)\|^2_{L^2}
   &\leq& \eps \|D u\|^2_{\bL(T)}
  +\varepsilon \E\sup_{t\leq T}\|u(t)\|^2_{L^2}+c(\varepsilon)\sum_{i=0}^d\|f^i\|^2_{\bL(T)} \\
  &\leq&  \eps \|D u\|^2_{\bL(T)}
  +\varepsilon \E\sup_{t\leq T}\|u(t)\|^2_{L^2}
  +c(\varepsilon)\|f\|^2_{\bH^{-1}(T)}. \end{eqnarray*}
 By Davis (first) inequality and L\'evy system,
 \begin{eqnarray}
\E \left[ \sup_{0\leq s\leq t}|M_s^k| \right] &\leq &
  2\sqrt{6} \, \E \left[ [M^k,M^k]_t^{1/2} \right]
 \leq    2\sqrt{6} \, \E \left[ \sum_{0<s\leq t}
\|g^k(t)\Delta Z^k_t\|^2_{L^2} \right] \nonumber \\
&\leq &  2\sqrt{6} \,  \wh {c}_{k}^2  \, \E \left[ \int_0^t
\|g^k(s)\|^2_{L^2} ds \right] \label{e:3.7}
\end{eqnarray}
and
\begin{eqnarray*}
&& \E \left[
 \sup_{0\leq s\leq t} \sum_{k=1}^{\infty} \Big| \int_0^s
(u(r-),g^k(r)) dZ^k_r \Big| \right]
\\
&\leq&  2 \sqrt{6} \, \sum_{k=1}^{\infty} \, \E \left[\left(
\sum_{0<s\leq t}
(u(s-),g^k(s))^2 (\Delta Z^k_s)^2\right) ^{1/2} \right] \\
&\leq& 2\sqrt{6} \,\sum_{k=1}^{\infty}\, \E \left[\sup_{s\leq
t}\|u(s)\|_{L^2}\left(\sum_{0<s\leq t} \|g^k(s)\|^2_{L^2} (\Delta
Z^k_s)^2\right) ^{1/2}\right] \\
&\leq&  \varepsilon \E\left[\sup_{s\leq
t}\|u(s)\|^2_{L^2}\right]+c(\varepsilon) \sum_{k=1}^{\infty}\E
\left[ \sum_{0<s\leq t}
\|g^k(s)\|^2_{L^2} (\Delta Z^k_s)^2  \right]\\
&\leq & \varepsilon \E\left[\sup_{s\leq
t}\|u(s)\|^2_{L^2}\right]+c(\varepsilon)\sum_{k=1}^{\infty} \wh
{c}_{k}^2 \, \E\int^t_0\|g^k (s)\|^2_{L^2}ds.
\end{eqnarray*}
It follows from \eqref{e:2.5} that
\begin{eqnarray*}
 \E \left[ \sup_{t\leq T} \| u(t)\|_{L^2}^2 \right]
 &\leq & \varepsilon \E\left[\sup_{s\leq T}\|u(s)\|^2_{L^2}\right]
 + \E\| u_0\|_{L^2}^2 +
 %\| Du\|^2_{\bL(T)}
  \eps \| Du\|^2_{\bL(T)}
 +c(\varepsilon) \| f\|^2_{\bH^{-1}(T)}
 +c(\varepsilon)\|g\|^2_{\bL(T,\ell^2)}.
 \end{eqnarray*}
Thus \eqref{e:2.4} is proved if one chooses  $\varepsilon\leq 1/2$.
Now we prove the completeness of the space $\cH^1(T)$. Let $\{u_n:
n=1,2,...\}$ be a Cauchy sequence in $\cH^1(T)$. Let $f_n:=\bD u_n$,
$g_n:=\bS u_n$ and $u_{n0}:=u_n(0)$. Then there exist $u\in
\bH^1(T)$, $f\in \bH^{-1}(T)$, $g\in \bL(T,\ell^2)$ and $u_0\in U^2$
so that $u_n$, $f_n$,
 $g_n=\{g_n^k, k\geq 1 \}$
and $u_{n0}$ converge to $u,f,g$ and $u_0$, respectively. Let
%$\phi\in C^{\infty}_0$ be fixed. By taking the limit from
 $\phi\in C^{\infty}_c$  be fixed. Since
$$
(u_n(t),\phi)=(u_{n0},\phi)+\int^t_0(f_n(s),\phi)ds+ \sum_{k\geq 1}
  \int^t_0(g^k_n (s),\phi) dZ^k_s,
$$
% we obtain for each $t$,
 taking $n\to \infty$, we have for each $t>0$,
\begin{equation}      \label{eqn 10.18-1}
(u(t),\phi)=(u_0,\phi)+\int^t_0(f(s),\phi)ds+ \sum_{k\geq 1}
\int^t_0(g^k (s),\phi) d Z^k_s \quad \hbox{a.s.}
\end{equation}
 Since we already proved
$$
 \E \left[ \sup_{t\leq T}\|u_n-u_m\|^2_{L^2} \right] \leq
 c\|u_n-u_m\|^2_{\cH^1(T)},
$$
we conclude that $(u_n(t),\phi)$ is uniformly Cauchy in $t\in [0, T]$,
(\ref{eqn 10.18-1}) holds for all $t\leq T$ a.s., and $u$ is right
continuous having  left limits in $L^2$ a.s. Consequently $u\in
\cH^1(T)$. \qed

\begin{assumption}
                         \label{A2.3}
{\rm (i)}  The coefficients $a^{ij}, \bar{b}^i, b^i, c,\sigma^{ik}$
and $\mu^k$ are $\cP \otimes \cB(\bR^d)$-measurable functions.

{\rm (ii)}  For each $\omega,t,x, i,j$,
$$
  |a^{ij}|+  |\bar{b}^i|+|b^i| +|c|+ \Big(
 \sum_{k=1}^{\infty}(\beta^2_k+\wh{c}^2_k)(|\sigma^{ik}|^2+|\mu^k|^2)\Big)^{1/2}
\leq K.
$$
\end{assumption}

\begin{lemma}(A priori estimate)\label{L:priori}
Let Assumptions  \ref{A2.1} and \ref{A2.3}  hold.
%Then  any solution
 Then for every solution
$u\in \cH^1(T)$ of equation \eqref{5.15.0}, we have
\begin{equation}
                    \label{111}
%\|u\|_{\cH^1(t)}\leq
 \|u\|_{\cH^1(T)}\leq
ce^{cT}\left(\|f\|_{\bH^{-1}(T)}+\|g\|_{\bL(T,\ell^2)}+\|u_0\|_{U_2}\right),
\end{equation}
where $c=c(\delta,K)$.
\end{lemma}

\pf  We proceed as in the proof of Theorem \ref{T:2.7}. As
before, rewrite $f\in \bH^{-1}(T)$ as
$$
f=f_0+\sum_{i=1}^d \frac{\partial}{\partial x_i}  f_i
\quad\hbox{with } \  f^i\in \bL(T)
$$
%so that
 and
$$
\sum_{i=0}^d\|f_i\|_{\bL(T)}\leq c\|f\|_{\bH^{-1}(T)}.
$$
As in the proof of Theorem \ref{T:2.7},
%we may
without loss of generality, we may and do
assume that $u,f,g$
are sufficiently smooth in $x$.  By $h^k$ we denote the predictable
version of $\sigma^{ik}u_{x^i}+\nu^ku+g^k$. By Ito's formula (cf.
\cite{HWY}), we have
\begin{eqnarray} \label{e:2.10}
  \E \left[\|u(t)\|^2_{L^2}\right] &=&
\E \left[\|u_0\|^2_{L^2} \right] +2\E \left[\int^t_0
\left(-(a^{ij}u_{x^j}+\bar{b}^iu+f_i,u_{x^i})_{L^2}+(b^iu_{x^i}+cu+f_0,u)_{L^2}
\right)ds \right] \nonumber\\
&& +\sum_k\beta^2_k\int^t_0\|h^k\|^2_{L^2}\,ds + \, 2\E
\left[\sum_k\int^t_0(h^k,u(s-))_{L^2}dZ^k_s
\right] \nonumber \\
&& + \, \sum_k\E \left[ \sum_{0<s\leq t}\|h^k\Delta Z^k_s\|^2_{L^2}
\right].
\end{eqnarray}
It is easy to show
\begin{eqnarray*}
%&&\sum_k\beta^2_k\int^t_0\|h^k\|^2_{L^2}\,ds=\E\sum_k
%\beta^2_k\int^t_0(\sum_i \sigma^{ik}u_{x^i}+\nu^ku+g^k)^2_{L^2}dt\\
 \E \left[\beta^2_k\int^t_0\|h^k\|^2_{L^2}\,ds\right]
 &=& \E\left[  \beta^2_k\int^t_0 \ \sigma^{ik}u_{x^i}+\nu^ku+g^k\|^2_{L^2}dt
  \right]\\
&\leq&
2\E \left[\int^t_0(\alpha^{ij}_1u_{x^i},u_{x^j})_{L^2}ds\right]
+\varepsilon\|Du\|^2_{\bL(t)}+c(\varepsilon)\|u\|^2_{\bL(t)}
+c(\varepsilon)\|g\|^2_{\bL(t,\ell^2)},
\end{eqnarray*}
where $\alpha^{ij}_1=\frac{1}{2}\sum_k
\beta^2_k\sigma^{ik}\sigma^{jk}$.  Also,
\begin{eqnarray*}
&&   \sum_k \E \left[\sum_{0<s\leq
t}\|h^k\Delta Z^k_s\|^2_{L^2} \right] \\
& =& \sum_k \wh {c}_k^2
\E \left[ \int^t_0\|\sigma^{ik}u_{x^i}+\mu^ku+g^k\|^2_{L^2}ds \right] \\
%removed:
%&=& \sum_k \wh{c}^2_k \E \Big[ \int^t_0 \Big(
%(\sigma^{ik}\sigma^{jk}u_{x^i},u_{x^j})+\|\wh{c}_k\mu^ku\|^2_{L^2}
%+\|\wh{c}_kg^k\|^2_{L^2}+(2\wh{c}_k\sigma^{ik}u_{x^i},\wh{c}_k\mu^ku)
% + 2(\wh{c}_k\sigma^{ik}u_{x^i},\wh{c}_kg^k) \\
% && \hskip 0.5truein +2(\wh{c}_k\mu^ku,\wh{c}_kg^k) \Big) ds  \Big]  \\
&\leq& 2 \E \left[ \int^t_0 ( \alpha^{ij}_2u_{x^i}, \, u_{x^j})ds
\right] +\varepsilon
\|Du\|^2_{\bL(t)}+c(\varepsilon)\|u\|^2_{\bL(t)}+c(\varepsilon)\|g\|^2_{\bL(t,\ell_2)},
\end{eqnarray*}
where
$\alpha^{ij}_2=\frac{1}{2}\sum_k\wh{c}^2_k\sigma^{ik}\sigma^{jk}$.
 Similarly,
%$$
%\E \left[ \int^t_0 \left(
%(\bar{b}^iu+f^i,u_{x^i})+(b^iu_{x^i}+cu+f^0,u)\right)ds \right] \leq
%\varepsilon \|Du\|^2_{\bL(t)}+c\|u\|^2_{\bL(t)}+c\|f^i\|^2_{\bL(t)}.
 \begin{eqnarray*}
&&\E \left[ \int^t_0 \left( (\bar{b}^i , u_{x^i})+\sum_{i=1}^d ( +f_i,u_{x^i}) +(b^iu_{x^i}+cu+f_0,u)\right)ds \right] \\
&\leq &
\varepsilon \|Du\|^2_{\bL(t)}+c(\eps) \|u\|^2_{\bL(t)}+c(\eps) \sum_{i=0}^d \|f_i\|^2_{\bL(t)}.
 \end{eqnarray*}
 Thus we have from \eqref{e:2.10} that for each $t\leq T$,
\begin{eqnarray*}
&& \E \left[\|u(t)\|^2_{L^2} \right] +2\E \left[ \sum_{i,j=1}^d
\int^t_0((a^{ij}-\alpha^{ij})u_{x^i}, \, u_{x^j}) \right]
ds \\
& \leq& \E \left[ \|u_0\|^2_{L^2} \right] +\varepsilon
\|Du\|^2_{\bL(t)}
%+c\int^t_0\E\|u(s)\|^2_{L^2}ds+c\|f^i\|^2_{\bL(t)}+c\|g\|^2_{\bL(t,\ell^2)}.
 +c(\eps) \int^t_0\E \left[\|u(s)\|^2_{L^2}\right] ds+c(\eps) \sum_{i=0}^d \|f_i\|^2_{\bL(t)}+c(\eps) \|g\|^2_{\bL(t,\ell^2)}.
\end{eqnarray*}
On the other hand, we know from condition \eqref{e:2.2} that
$$   \sum_{i,j=1}^d
 ((a^{ij}-\alpha^{ij})u_{x^i}, \, u_{x^j}) \geq \delta \| D u \|_{L^2}^2.
$$
 The above two displays
 % yield
 together with Grownwell's inequalty yield
$$
\|u\|_{\bH^1(T)}\leq ce^{cT}\left(\|u_0\|_{U_2}+
\|f\|_{\bH^{-1}(T)}+\|g\|_{\bL(T,\ell^2)}\right),
$$
where $c=d(\delta,K)$. The lemma is proved. \qed

\begin{remark}
                \label{neq remark}
The proof of Lemma \ref{L:priori} shows that if
$\bar{b}^i=b^i=c=\nu^k=0$, then
$$
\|u_x\|_{\bL(T)}\leq
c\left(\|f\|_{\bH^{-1}(T)}+\|g\|_{\bL(T,\ell^2)}+\|u_0\|_{U_2}\right)
$$
where $c$ is {\bf{independent}} of $T$.
\end{remark}

 \vspace{5mm}

For $\lambda \in [0,1]$, denote
$$
a^{ij}_{\lambda}=\lambda a^{ij}+(1-\lambda)\delta^{ij}, \quad
\sigma^{ik}_{\lambda}=\lambda \sigma^{ik},
$$
$$
\bar{b}^i_{\lambda}=\lambda \bar{b}^i, \quad b^i_{\lambda}=\lambda
b^i,\quad  c_{\lambda}=\lambda c, \quad \mu^k_{\lambda}=\lambda
\mu^k.
$$
$$
L_{\lambda}u:=\lambda Lu+(1-\lambda) \Delta u=
\frac{\partial}{\partial
x_i}(a^{ij}_{\lambda}u_{x^j}+\bar{b}^i_{\lambda})+b^i_{\lambda}u_{x^i}+c_{\lambda}u,
$$
$$
\Lambda^k_{\lambda}u:=\lambda \Lambda^k u
:=\sigma^{ik}_{\lambda}u_{x^i}+\mu^k_{\lambda}u \qquad \hbox{for }
k\geq 1.
$$
Note that
$$
L_{\lambda_1} u-L_{\lambda_2}u=(\lambda_1-\lambda_2)(L-\Delta)u,
\qquad
\Lambda_{\lambda_1}u-\Lambda_{\lambda_2}u=(\lambda_1-\lambda_2)\Lambda
u,
$$
where $\Lambda_\lambda u :=(\Lambda^1_\lambda u,  \Lambda^2_\lambda
u, \cdots )$,   $\Lambda  u :=(\Lambda^1  u, \Lambda^2  u, \cdots
)$, and
\begin{equation}
                                 \label{eqn 11.06}
\|L_{\lambda_1} u- L_{\lambda_2}u\|_{H^{-1}}+
\|\Lambda_{\lambda_1}u-\Lambda_{\lambda_2}u\|_{L^2(\ell^2)}\leq
c|\lambda_1-\lambda_2|\|u\|_{H^1}.
\end{equation}

\begin{remark}  \label{remark 5} \rm
It is trivial to check that  a priori estimate (\ref{111}) holds
with the same constant $C$ if $u\in \cH^1(T)$ is a solution of the
equation obtained by replacing the  coefficients
$a^{ij},\bar{b}^i,\cdots, \mu^k$  in (\ref{5.15.0}) by
$a^{ij}_{\lambda},\bar{b}^i_{\lambda},\cdots, \mu^k_{\lambda}$,
respectively, for every $\lambda \in [0, 1]$.
\end{remark}

Here is the main result of this section.

\begin{thm}  \label{main thm}
 Suppose {\bf Assumptions \ref{A2.1}} and {\bf \ref{A2.3}} hold.
Then  for every  $f\in \bH^{-1}(T)$, $g\in \bL(T, \ell^2)$ and
$u_0\in U_2$, equation (\ref{5.15.0}) has a unique solution $u\in
\cH^1(T)$
 with $u(0)=u_0$,
and
\begin{equation}
                         \label{e: main}
\|u\|_{\cH^1(T)}\leq ce^{cT}(\|f\|_{\bH^{-1}(T)}+\|g\|_{\bL(T,
\ell^2)}+\|u_0\|_{U_2}),
\end{equation}
where $c=c(\delta,K)$.
\end{thm}

 \pf In view of the a priori estimate in Lemma \ref{L:priori}, it
suffices to show that there is a  solution
 to \eqref{5.15.0}.
First, we show that for any given $f\in \bH^{-1}(T), g\in
\bL(T,\ell^2)$ and $u_0\in U_2$, the equation
\begin{equation}
                          \label{eqn 11.02}
du=(\Delta u+f)dt + g^k dZ^k_t, \quad u(0)=u_0
\end{equation}
has a solution $u\in \cH^1(T)$. Due to a priori estimate
(\ref{111}), Remark \ref{remark 5}
 and standard approximation argument, we may assume that
$f,u_0$ are infinitely differentiable in $x$ with compact
supports. Also by the same reasoning (also see Theorem 3.10 in
\cite{Kr99}), we may assume that  $g^k=0$ for all $k\geq N$
 for some $N\geq 1$,
 and $g^k$ is of the type
$$
g^k(t)=\sum_{i=1}^m I_{(\tau_{i-1},\tau_i]}(t)
 \varphi_i (x),
$$
where $\tau_i$ are bounded stopping times and
 $\varphi_i \in C^\infty_c (\bR^d)$. Define
$$
v(t)=\sum_{k=1}^N \int^t_0 g^k (s) dZ^k_s .
$$
Then it is easy to
 see that
$v\in \cH^1(T)$. Note that $u$ satisfies
(\ref{eqn 11.02}) if and only if $\bar{u}=u-v$ satisfies
$$
d\bar{u}=(\Delta \bar{u} +\Delta v+f)dt  \qquad \hbox{with} \quad  \bar{u}(0)=u_0.
$$
Since this equation has a solution in $\cH^1(T)$ (see Theorem 5.1 in
\cite{Kr99}), we conclude that equation (\ref{eqn 11.02}) has a
solution $u$ in $\cH^1(T)$.

Let $J\subset [0,1]$ denote the set of $\lambda$, so that for any
$f,g,u_0$, the equation
\begin{equation}
                            \label{11}
du=(L_{\lambda}u+f)dt+(\Lambda^k_{\lambda} u+g^k)dZ^k_t, \quad
u(0)=u_0
\end{equation}
has a solution $u\in \cH^1(T)$.  Then as proved above,
 $0\in J$.
Now assume $\lambda_0\in J$, and note that $u$ is a solution
of equation (\ref{11}) if and only if
\begin{equation}\label{e:2.12}
du=(L_{\lambda_0}u+(L_{\lambda}u-L_{\lambda_0}u+f))dt+(\Lambda_{\lambda_0}u+
(\Lambda^k_{\lambda} u-\Lambda^k_{\lambda_0}u+g^k))dZ^k_t.
\end{equation}
 Note
that   $D:H^n \to H^{n-1}$ is a bounded operator. Thus  for
any $u\in \cH^1(T)$, $k\geq 1$ and $\lambda\in [0, 1]$, we have
$$
L_{\lambda}u\in \bH^{-1}(T)  \qquad  \hbox{and} \qquad
\Lambda_{\lambda}u \in \bL(T,\ell^2).
$$
Recall $\lambda_0\in J$. Denote  $u^0=u_0$ and for $n=1,2,\cdots$ we
define $u^{n+1}\in \cH^1(T)$  as the solution of the equation
$$
du^{n+1}=(L_{\lambda_0}u^{n+1}+f_n)dt +(\Lambda_{\lambda_0}u^{n+1} +
g^k_n )dZ^k_t, \quad u^{n+1}(0)=u_0
$$
where
$$
f_n:=L_{\lambda}u^n-L_{\lambda_0}u^n+f \quad \hbox{and} \quad
g^k_n:=\Lambda^k_{\lambda} u^n-\Lambda^k_{\lambda_0}u^n+g^k.
$$
 By Remark \ref{remark 5} and  inequality (\ref{eqn 11.06}), we have
 \begin{eqnarray*}
\|u^{n+1}-u^n\|_{\cH^1(T)} &\leq&
c\|(L_{\lambda}-L_{\lambda_0})(u^n-u^{n-1})\|_{\bH^{-1}(T)}
+c\|(\Lambda_{\lambda}-\Lambda_{\lambda_0})(u^n-u^{n-1})\|_{\bL(T)}
\\
& \leq& c \|\lambda-\lambda_0|\|u^n-u^{n-1}\|_{\bH^1(T)}.
\end{eqnarray*}
Let $\eps_0=c/2$.  Then for $\lambda \in (\lambda_0 -\eps_0, \,
\lambda+\eps_0)$, $\|u^{n+1}-u^n\|_{\cH^1(T)}\leq \frac12
|\lambda-\lambda_0|\|u^n-u^{n-1}\|_{\bH^1(T)}$ for every $n\geq 1$
and so  $u^n$ converges to some $u$ in $\cH^1(T)$.  It follows that
$u$ solves equation \eqref{e:2.12}. This proves that
$(\lambda_0-\eps_0,\, \lambda_0+\eps_0)\cap [0, 1]\subset  J$.
Consequently we conclude $J=[0,1]$. \qed

The following remark plays the key role when we weaken condition (\ref{eqn e.1}) later in the next section.

\begin{remark}      \label{remark 12.10}
Let $\tau\leq T$ be a stopping time.
 We use $1_{[\![0, \tau]\!]}$ to denote the random
 process $t\mapsto 1_{[0, \tau]}(t)$.
 For an $H^1$-valued
$\cP^{dP\times dt}$-measurable process $u$, write $u\in \bH^1(\tau)$
if
$$\|u\|^2_{\bH^1(\tau)}:=\E \left[\int^{\tau}_0 \|u\|^2_{H^1}ds\right]
<\infty.
$$
Similarly define $\bL(\tau,\ell_2)$ and $\cH^1(\tau)$. Then Theorem
\ref{main thm}  holds
 with the deterministic time $T$ replaced by the stopping time $\tau$.
 Indeed, the existence of solution $u\in \cH^1(\tau)$ and the estimate
(\ref{e: main})  are easily obtained by applying
%the theorem
 Theorem \ref{main thm}
with $\bar{f}=f1_{[\![ 0, \tau]\!]}$ and $\bar{g}=g 1_{[\![0, \tau]\!]}$
  in place of $f$ and $g$,
respectively. Now let $u\in \cH^1(\tau)$ be a solution. According to
Theorem \ref{main thm} we can define $v\in \cH^1(T)$ as the solution
of
\begin{equation}
                 \label{e: remark 12.10}
  dv=(\Delta v+(\bD u-\Delta u)1_{[\![ 0, \tau]\!]})dt+ 1_{[\![0, \tau ]\!]} \bS^ku \,dZ^k_t, \quad v(0)=u(0).
\end{equation}
Then for $t\leq [0,\tau)$, $d(u-v)=\Delta(u-v)dt$ and therefore we
conclude that $u(t)=v(t)$ for all $t\leq \tau $, a.s.. Thus,
equation (\ref{e: remark 12.10}) becomes
\begin{eqnarray} \label{e: remark 12.10.1}
dv&=&\left( \sum_{i =1}^d \frac{\partial}{\partial x^i} \left( \sum_{j=1}^d a^{ij}_{\tau}v_{x^j}
 + \bar{b}^i_{\tau} v \right)+b^i_{\tau}v_{x^i}+c_{\tau}v+f
  1_{[\![ 0, \tau]\!]} \right)dt \nonumber \\
&& +\sum_{k\geq 1} \left( \sum_{i=1}^d \sigma^{ik}_{\tau}v_{x^i}+\mu^k_{\tau}v+g^k
  1_{[\![ 0, \tau]\!]} \right)
 dZ^k_t,
\end{eqnarray}
where
$$
a^{ij}_{\tau}=a^{ij}
 1_{[\![ 0,  \tau ]\!]}+\delta^{ij} 1_{ ]\!] \tau, \infty [\![}, \quad
\bar{b}^i_{\tau}=\bar{b}^i
 1_{[\![ 0, \tau]\!]},\quad b^i_{\tau}=b^i1_{[\![ 0, \tau]\!]}, \quad \cdots,
 \mu^k_{\tau}=\mu^k1_{[\![ 0, \tau]\!]}.
$$
The uniqueness result of equation \eqref{5.15.0} in the space
$\cH^1(\tau)$ follows from  the uniqueness result of equation
\eqref{e: remark 12.10.1} in $\cH^1(T)$.
\end{remark}

\mysection{Some extensions}
                                           \label{section extension}
In this section we give two extensions of Theorem  \ref{main thm}.
First, we consider the nonlinear equation
\begin{eqnarray}    \label{e:3.1}
du &=& \left(\sum_{i=1}^d \frac{\partial}{\partial x_i} \left(
\sum_{j=1}^d a^{ij}u_{x^j}+
\bar{b}^iu \right)+ \sum_{i=1}^d b^iu_{x^i}+ cu+ f(u) \right)\,dt
\nonumber\\
&&  + \sum_{k\geq 1}
\left(\sum_{i=1}^d \sigma^{ik} u_{x^i}+\mu^{k}u+g^{k}(u)\right) dZ^{k}_{t},
\end{eqnarray}
where $f(u)=f(\omega,u,t,x)$ and $g^k(u)=g^k(\omega,u,t,x)$.

\begin{assumption} \label{A3.1} \rm \begin{description}
 \item{(i)} For any $u\in \bH^1(T)$,
 $$
 f(u)\in \bH^{-1}(T) \quad \hbox{and} \quad
  g(u):= (g^1(u), g^2(u), \cdots )
 \in \bL(T,\ell^2).
 $$

 \item{(ii)} For  every $\varepsilon >0$, there exists a constant
 $K_1=K_1(\varepsilon)$  so that for every $t\in (0, T]$ and  $u,v\in \bH^1(t)$,
\begin{equation}\label{e:3.2}
 \|f(u)-f(v)\|^2_{\bH^{-1}(t)}+\|g(u)-g(v)\|^2_{\bL(t,\ell^2)}\leq
 \varepsilon \|u-v\|^2_{\bH^1(t)}+K_1\|u-v\|^2_{\bL(t)}.
 \end{equation}
 \end{description}
\end{assumption}

\begin{thm}
                      \label{main thm-2}
 Suppose {\bf Assumptions \ref{A2.1}}, {\bf \ref{A2.3}} and
 {\bf \ref{A3.1}}
 hold. Then for any $u_0\in U_2$,
equation  \eqref{e:3.1} with initial data $u_0\in U_2$ has a unique
solution $u\in \cH^1(T)$, and
\begin{equation}\label{e:3.3}
\|u\|_{\cH^1(T)}\leq c(\|f(0)\|_{\bH^{-1}(T)}+\|g(0)\|_{\bL(T,
\ell^2)}+\|u_0\|_{U_2})
\end{equation}
where $f(0)=f(\omega,0,t,x), g(0)=g(\omega,0,t,x)$ and
$c=c(\delta,K,T)>0$.
\end{thm}

\pf
 We will use fixed point theorem to show the existence and
uniqueness of the solution to \eqref{e:3.1}. Estimate \eqref{e:3.3}
follows from \eqref{111}, condition \eqref{e:3.2} and the
Grownwall's inequality. Let $\cR(f,g)\in \cH^1(T)$ denote the
solution of (\ref{5.15.0}) with  initial data $u_0$. Then by Theorem
\ref{main thm},
$$
 \cR u:=\cR(f(u),g(u)) \qquad \hbox{for } u\in \cH^1(T)
$$
is   well defined and $\cR$ is a map from $\cH^1(T)$ to $\cH^1(T)$.
 Define $u^0=\cR(f(0), g(0))$
and $u^{n+1}=\cR(f(u^n),g(u^n))$. Then
by Theorem \ref{main thm} and assumption \eqref{e:3.2}, for any $t\leq T$,
\begin{eqnarray*}
\|\cR u-\cR v\|^2_{\cH^1(t)} &\leq&
c\varepsilon\|u-v\|^2_{\cH^1(t)}+cK_1\|u-v\|^2_{\bL(t)} \\
& \leq& c\varepsilon\|u-v\|^2_{\cH^1(t)}+cK_1 \int_0^t
\|u-v\|^2_{\cH^1(s)}ds
\end{eqnarray*}
where the last inequality is from \eqref{e:2.4}. Taking
$\eps=1/(2c)$ and then letting $t_0>0$ be small enough so that
$cK_1t_0<1/4$,  we have
 \begin{equation}\label{e:3.4}
  \|\cR u-\cR v\|^2_{\cH^1(t_0)} \leq \frac12 \|u-v\|^2_{\cH^1(t_0)}
\end{equation}
This contraction implies that $u^n$ converges to some $u$ in
$ \cH^1(t_0)$ and $u$ is a solution to \eqref{e:3.1} on $[0, t_0]$
 with initial value $u_0$.
The inequality \eqref{e:3.4} further implies
the uniqueness of solution to \eqref{e:3.1} on $[0, t_0]$
  initial value $u_0$.
 Iterating this procedure at most $[T/t_0]+1$ many time intervals
of size no larger than $t_0$ and using estimate \eqref{e:2.4},
we get the desired results on time interval $[0, T]$.
  For more details, we refer the reader to the proof of Theorem 6.4 in \cite{Kr99}. \qed

\begin{example}
Let's consider an equation with fractional Laplacian. For simplicity
assume $g^k(u)=0$ for $k\geq 2$. Take $f(u)=(-\Delta)^{\alpha/2}u$
and $g(u)=g^1(u)=(-\Delta)^{\beta/2}u$ where $\alpha<2$ and
$\beta<1$, then obviously for any $\varepsilon>0$,
$$
\|f(u)-f(v)\|^2_{\bH^{-1}(t)}+\|g(u)-g(v)\|^2_{\bL(t)}\leq
c\|u-v\|^2_{\bH^{-1+\alpha}(t)}+c\|u-v\|^2_{\bH^{\beta}(t)}
$$
$$
\leq \varepsilon \|u-v\|^2_{\bH^1(t)}+K_1\|u-v\|^2_{\bL(t)},
$$
where for the second inequality we use the following fact: if
$\gamma=\kappa \gamma_1+(1-\kappa)\gamma_0$ and $\kappa\in [0,1]$
then $\|u\|_{H^{\gamma}}\leq
N\|u\|^{\kappa}_{H^{\gamma_1}}\|u\|^{1-\kappa}_{H^{\gamma_0}}$.
Thus the existence and uniqueness of equation (\ref{e:3.3}) in $\cH^1(T)$ is guaranteed by Theorem \ref{main thm-2}.

\end{example}

For a stopping time $\tau \in (0,T]$ write $u\in \bH^1_{\loc}(\tau)$ if there exists a
sequence of stopping times $\tau_n \uparrow \infty$ so that $u\in
\bH^1(\tau \wedge \tau_n)$ for each $n$.

The following is a weakened version of {\bf Assumption \ref{A2.1}}.
 \begin{assumption}
            \label{A3.3}
 There exists an integer $N_0\geq 1$ so that

 {\rm (i)}  $\wh{c}_k<\infty$ for all
 integer $k>N_0$;

{\rm (ii)}  for some $\delta >0$,
\begin{equation}
                                \label{eqn 7.23}
( a^{ij}-\alpha^{ij}_{N_0} )_{d\times d}>\delta I_{d\times d},
\end{equation}
where $\alpha^{ij}_{N_0}:= \frac{1}{2} \sum_{k=N_0+1}^\infty (\wh
{c}_{k}^2+\beta^2_k) \sigma^{ik}\sigma^{jk}$.
\end{assumption}

Here is our second
extension.

\begin{thm}
Let Assumption \ref{A3.3} hold and $\sigma^i_k=0$ for $k\leq N_0$.
 Then for any
  $u_0\in U_2$,
$f\in \bH^{-1}(T)$
%.
 and  process $g=(g^1,g^2,\cdots)$ having entries in $\bL_2(T)$ so that $\sum_{N_0 +1}^{\infty}\wh{c}^2_k\|g^k\|^2_{\bL(T)}<\infty$,  there exists unique $u\in \bH^1_{\loc}(T)$
such that

{\rm (i)} $u(t)$ is right continuous with left limits in $L^2$ $a.s.$,

{\rm (ii)} for any $\phi\in C^\infty_c (\bR^d)$, the equality
\begin{eqnarray}
(u(t),\phi) &=&(u_0,\phi)+\int^t_0 \left((-a^{ij}u_{x^j}-\bar{b}^iu,
\,
\phi_{x^i})+(b^iu_{x^i}+cu+f, \, \phi)\right)\,ds \nonumber \\
&&
+\int^t_0\left((\sigma^{ik}u_{x^i},\phi)+
 (\mu^k,\phi) + (g^k, \phi) \right)\,dZ^k_s
             \label{e:3.5}
\end{eqnarray}
holds for all $t<T$ $a.s.$.
\end{thm}

\pf {\bf{Step 1}}.
 First assume that  Assumption \ref{A2.1} holds, that is, $\wh c_k<\infty$ for each $k$.  Let $\tau \leq T$ be a stopping time. We show
that the pathwise solution is unique in $\bH^1_{\loc}(\tau)$.
Let $u\in \bH^1_{\loc}(\tau)$ be a path-wise solution, that is, $u$
satisfies the conditions (i) and (ii) in the theorem for $t < \tau$. Define
$\tau_n=\tau \wedge \inf\{t:\int^t_0\|u\|^2_{H^1}ds>n\}$.  Then $u\in
\bH^1(\tau_n)$ and $\tau_n \uparrow \tau$ since
$\int^t_0\|u\|^2_{H^1}ds<\infty$ for all $t <\tau$, a.s.  By
Remark \ref{remark 12.10},
$$
\|u\|_{\bH^1(\tau_n)}\leq
c(T,d,K)(\|f\|_{\bH^{-1}(\tau_n)}+\|g\|_{\bL(\tau_n,\ell_2)}+\|u(0)\|_{U_2}).
$$
By letting $n \to \infty$ we find that $u\in \cH^1(\tau)$, and the
uniqueness of the pathwise solution under  Assumption \ref{A2.1}
follows from Remark \ref{remark 12.10}. Note that the existence of pathwise solution under Assumption \ref{A2.1} in $\bH^1(\tau)$ also follows from Theorem \ref{main thm}.

\medskip

{\bf{Step 2}}. For the general case,  note that for each $n>0$
$$ \wh c_{k,n}:= \left(\int_{\{z\in \bR: |z|\leq n\}} |z|^2 \nu_k (dz) \right)^{1/2} \qquad \hbox{for } k\leq N_0.
$$
Consider the  L\'evy
processes $(Z^1_n, \cdots,  Z^{N_0}_n, Z^{N_0+1},\cdots)$ in place of $(Z^1, Z^2\cdots)$, where $Z^k_n (k\leq N_0)$ is a L\'evy
process obtained from $Z^k$ by removing all the jumps that has
absolute size strictly large than $n$. Note that condition (\ref{e:2.2}) is valid with $\wh {c}_k$  replaced by $\wh {c}_{k,n}$  since $\sigma^{ik}$ are
assumed to be zero for $k\leq N_0$.

By Step 1,
 there is a
unique pathwise solution $v_n\in  \cH^1(T)$
 with  $Z^k_n$ in place of $Z^k$ for
$k=1, 2, \cdots, N_0$. Let $T_n$ be the first time that one of the
L\'evy processes $\{Z^k, 1\leq k\leq N_0\}$ has a jump of (absolute)
size in $(n, \infty)$. Define $u(t)=v_n(t)$ for $t <T_n\wedge T$. Note that for $n<m$, by Step 1, we have $v_n(t)=v_m(t)$  for $t< T_n$. This is because, for $t<T_n$,  both $v_n$ and $v_m$
satisfies (\ref{e:3.5}) with each term inside the stochastic integral multiplied by $1_{s<T_n}$ (and with $Z^k_n$, $k\leq N_0$, in place of $Z^k$).
Thus $u$ is well defined. 
By letting $n\to \infty$, one constructs  unique  pathwise solution $u$ in $\bH^1_{loc}(T)$. The theorem is proved.

\qed

\end{document}